\documentclass[12pt, reqno]{amsart}

\usepackage{amsmath}

\usepackage{palatino}

\usepackage{amssymb}

\usepackage{textcomp}

 2

\theoremstyle{remark}

\newcommand{\Q}{{\mathbb Q}}

\newcommand{\C}{{\mathbb C}}

\newcommand{\F}{{\mathbb F}}

\begin{document}

\title[Non-vanishing of Dirichlet series]{Non-vanishing of Dirichlet series with periodic coefficients}

\author[Tapas Chatterjee and M. Ram Murty]{Tapas Chatterjee\textsuperscript{1} and M. Ram Murty\textsuperscript{2}}

\address[Tapas Chatterjee and M. Ram Murty]
      {Department of Mathematics and Statistics,
       Queen's University, Kingston, Ontario,
       Canada, K7L3N6.}
\email[Tapas Chatterjee]{tapasc@mast.queensu.ca}

\email[M. Ram Murty]{murty@mast.queensu.ca}



\maketitle

\footnotetext[1]{Research of the first author was supported by a postdoctoral fellowship at Queen's University.}

\footnotetext[2]{Research of the second author was supported by an NSERC Discovery grant.}

\begin{abstract}

For any periodic function $f:{\mathbb N} \to {\mathbb C}$ with period $q$, we 
study the Dirichlet series $L(s,f):=\sum_{n\geq 1} f(n)/n^s.$  
It is well-known that this admits an analytic continuation
to the entire complex plane except at $s=1$,
where it has a simple pole with residue 
$$\rho:= q^{-1}\sum_{1\leq a\leq q} f(a).$$ Thus,
the function is analytic at $s=1$ when $\rho=0$ and in this case, 
we study its non-vanishing using the theory
of linear forms in logarithms and  Dirichlet $L$-series.
In this way, we give new proofs of an old criterion of Okada for the non-vanishing of $L(1,f)$
as well as a classical theorem of Baker, Birch and Wirsing. We also give some new necessary and sufficient 
conditions for the non-vanishing of $L(1,f)$.
\end{abstract}

\section{\bf Introduction}

\bigskip

Let $q$ be a positive integer and $f$ be a complex-valued periodic function with period $q$ 
which is not identically zero. The Dirichlet $L$-function $L(s,f)$ associated with 
$f$ is defined by the series

\begin{equation}\label{1}
L(s,f):=\sum_{n=1}^\infty \frac{f(n)}{n^s}.
\end{equation}

Since $f$ is periodic with period $q$, the above series can be written as

\begin{equation*}
L(s,f)=q^{-s}\sum_{a=1}^q f(a)\zeta(s,a/q), ~{\rm for \, } ~\Re(s)>1,
\end{equation*}

where $\zeta(s,x)$ is the Hurwitz zeta function defined for $\Re(s) >1$ by

$$\zeta(s,x) = \sum_{n=0}^\infty {1 \over (n+x)^s}, \qquad {\rm for }\quad 
 0 < x\leq 1. $$

Hurwitz \cite{hurwitz} showed that $\zeta(s,x)$ extends analytically to the entire complex plane except 
at $s=1$, where it has a simple pole with residue 1.

This shows that $L(s,f)$ extends analytically to the whole complex plane with a possible simple pole at
$s=1$ with residue $q^{-1}\sum_{\substack{a=1}}^{q} f(a)$. Hence $L(s,f)$ is an entire function if and only if
$\sum_{\substack{a=1}}^{q} f(a)=0$. From now on we will assume $\sum_{\substack{a=1}}^{q} f(a)=0$.  
In this case, we would like to know whether $L(1,f)=0$.

Partly motivated by prime number theory, this question was first raised
by Chowla \cite{chowla} in the case that $q$ is prime
and $f$ is rational-valued.  

If $f$ is algebraic-valued, its non-vanishing was established by Baker, Birch
and Wirsing \cite{BBW} under the additional conditions that $f(a)=0$
for each $a$ (mod $q$) satisfying $1 < (a,q) < q$ and the field
generated by the values of $f$ is disjoint from the $q$-th cyclotomic
field.  In this way, they answered the question of Chowla since the 
conditions are satisfied in the case $q$ is prime and $f$ is rational-valued.
Okada \cite{TO} was the first to derive necessary and sufficient
conditions on $f$ to ensure that $L(1,f)\neq 0$.  His derivation is long
and complicated.  In this paper, we offer another approach using
Dirichlet's $L$-functions and deduce Okada's criterion
from analytic properties of the classical $L$-functions.  We also
revisit the approach of Baker, Birch and Wirsing.  We are hopeful
that this approach will find new applications such as in the study
of the folklore Erd\"{o}s conjecture discussed in \cite{CM}. 

\smallskip

Since $f$ is a periodic function, we can define the Fourier transformation of $f$ by 

\begin{equation*}
\hat{f}(b)=\frac{1}{q}\sum_{a=1}^q f(a)\zeta^{-ab}_q 
\end{equation*}

where $\zeta_q=e^{\frac{2\pi i}{q}}$ and hence we have the Fourier inversion formula

\begin{equation*}
f(b)=\sum_{a=1}^q \hat{f}(a)\zeta^{ab}_q. 
\end{equation*}

It is also convenient to introduce the inner product on the
group of coprime residue classes (mod $q$):

$$(f,g) := {1 \over \varphi(q)} \sum_{\substack{a=1\\(a,q)=1}}^q f(a) \overline{g(a)}. $$

Henceforth, we suppose $f$ is rational-valued.  Then, we have $\hat{f}(b)\in \Q(\zeta_q)$.

Let $1=\omega_1,\cdots,\omega_{\varphi(q)}$ be an integral basis of $\Q(\zeta_q)$
over $\Q$. Then $\hat{f}(b)$ can be written as

\begin{equation}\label{2}
\hat{f}(b)=\sum_{j=1}^{\varphi(q)} c_{bj}\omega_j
\end{equation}
where $c_{bj}$'s are rational numbers. Note that, $\hat{f}(q)=0$ as $\sum_{a=1}^q f(a)=0$.

\medskip

Here is an outline of the results of this paper.
In section 3, we prove new necessary and sufficient conditions for the non-vanishing of
$L(1,f)$ in terms of the following theorems:

\smallskip
\noindent

\textbf{Theorem 1}. {\sl  $L(1,f)=0$ if and only if 

\begin{equation}\label{3}
\underset{{b=1}}{\prod^{q-1}}(1-\zeta^b_q)^{c_{bj}}=1
\end{equation}
for all $j=1,\cdots,\varphi(q)$.}

\smallskip
\noindent

\textbf{Theorem 2}.  {\sl $L(1,f)=0$ if and only if
$L(1,f_e)=0$ and $L(1,f_o)=0$, where $f_e$ and $f_o$ are the even and odd part of $f$ respectively.}

\medskip

Let $M(q)$ be the monoid generated by all the
prime divisors of $q$.  We prove:

\smallskip
\noindent

\textbf{Theorem 3.} {\sl  $L(1,f)=0$ if and only if 
\begin{equation*}
 \sum_{b\in M(q)}{f(ab)\over b}=0
\end{equation*}

for every $a$ with $1\le a <q, (a,q)=1$, and

\begin{equation*}
\sum_{b\in M(q)} \frac{(f_b,\chi_0)}{b}\log b=0, 
\end{equation*}

where $\chi_0$ is the principal Dirichlet character mod $q$ and $f_b(a)=f(ab)$.}

This theorem is reminiscent of Okada's criterion \cite{TO} and
in section 5, we establish the equivalence of Theorem 3 and Okada's criterion. In 
particular we prove the following theorem:

\smallskip
\noindent
\textbf{Theorem 4.}{\sl 
\begin{equation*}
\sum_{b\in M(q)} \frac{(f_b,\chi_0)}{b}\log b=0 
\end{equation*}

if and only if 

\begin{equation*}
 \sum_{\substack{r=1\\ (r,q)>1}}^qf(r)\epsilon(r,p)=0
\end{equation*}

for every prime divisor $p$ of $q$, where

$$
 \epsilon(r,p)=\begin{cases} 
v_p(r)  & \text{ if $v_p(r)<v_p(q)$,}\\

v_p(q)+\frac{1}{p-1}  & \text{ otherwise }
\end{cases}
$$

and for any integer $r$, $v_p(r)$ is the exponent of $p$ dividing $r$.}

\medskip

These questions and conjectures have a long history.
As the Dirichlet series \eqref{1} in general does not have an Euler product, even the 
existence of zeros in the domain of absolute convergence $\Re(s) >1$ cannot be ruled out. In 1969, as mentioned earlier, 
S. Chowla \cite{chowla} asked the question whether there exists a rational valued periodic function $f$ with prime 
period, such that $L(1,f)=0$. In 1973, Baker, Birch and Wirsing \cite{BBW} proved using Baker's theory 
of linear forms of logarithms, the following proposition which answered the question of Chowla.

\smallskip
\noindent
\textbf{Proposition 1}. {\sl If $f$ is a non-vanishing function defined on the integers with algebraic values and 
period $q$ such that (i) $f(r) = 0$ if $1 < (r, q) < q$,
(ii) the $q$-th cyclotomic polynomial $\Phi_q$ is irreducible over $\Q(f(1),\cdots, f (q))$, then 
\begin{equation*}
L(1,f)=\sum_{n=1}^\infty \frac{f(n)}{n}\ne 0. 
\end{equation*}}

In the last section, we give a new proof of the above proposition.

\par

One can also study these questions at points other than $s=1$, and this
was the focus of study in \cite{gmr}.  In addition, the authors 
of \cite{gmr} obtain
there a new generalization of the theorem of Baker, Birch and Wirsing over
number fields particularly in cases where the field generated by the 
values of $f$ is not necessarily disjoint from the $q$-th cyclotomic field.

\smallskip

\section{\bf Notations and Preliminaries}

\medskip 

In this section, we collect for the convenience of exposition, 
several results that will be used in the paper.
From now onwards, we denote the field of rationals by $\Q$, algebraic numbers by $\overline{\Q}$ and a 
number field by $\F$. The digamma function $\psi$ is the logarithmic derivative of the classical gamma 
function and is defined by the series

\begin{equation}\label{di}
 -\psi(x)=\gamma+{1\over x}+\sum_{n=1}^\infty\left({1\over {n+x}}-{1\over n}\right).
\end{equation}

\subsection{\bf A quick review of Baker's theory}

\bigskip

The following theorems due to Baker (see Theorem 2.1 and 2.4 of \cite{AB}) will play a crucial role in proving some 
of the theorems.

\smallskip
\noindent

\textbf{First version}. {\sl If $\alpha_1,\cdots, \alpha_n$ are non-zero algebraic numbers such that $\log\alpha_1, \cdots, \log\alpha_n$ 
are linearly independent over the field of rational numbers, then $1,\log\alpha_1, \cdots, \log\alpha_n$ are linearly independent 
over the field of algebraic numbers.}

\smallskip
\noindent

\textbf{Second version}. {\sl $\alpha_1^{\beta_1},\cdots,\alpha_n^{\beta_n}$ is transcendental for any algebraic numbers
$\alpha_1,\cdots,\alpha_n$ other than $0$ or $1$, and any algebraic numbers $\beta_1,\cdots,\beta_n$ with 
$1,\beta_1,\cdots,\beta_n$ linearly independent over the field of rationals. }

\bigskip

\subsection{\bf Okada's criterion}

 \medskip

In 1986, Okada \cite{TO} proved a proposition about 
the non-vanishing of $L(1,f)$ and 
Saradha and Tijdeman \cite{ST} modified his proposition, which we call Okada's criterion for the non-vanishing of
$L(1,f)$. Here is the proposition:

\smallskip
\noindent

\textbf{Proposition 2}. {\sl Let the $q$-th cyclotomic polynomial $\Phi_q$ be irreducible over $\Q(f(1),\cdots, f (q))$. 
As before, let $M(q)$ be the set of positive integers which are composed of prime factors of $q$.
Then $L(1,f)=0$ if and only if 
\begin{equation*}
 \sum_{m\in M(q)}{f(am)\over m}=0
\end{equation*}
for every $a$ with $1\le a <q, (a,q)=1$, and 

\begin{equation*}
 \sum_{\substack{r=1\\ (r,q)>1}}^qf(r)\epsilon(r,p)=0
\end{equation*}

for every prime divisor $p$ of $q$.}

We record here a variation of Lemma 5 in \cite{murty-murty} 
(see also \cite{MS}) which 
is useful in the last section of this paper, but is also of independent
interest.  

\smallskip

\textbf{Proposition 3}.  {\sl Let $\alpha_1, \alpha_2, ..., \alpha_n$ be positive
units in a number field of degree $>1$.  
Let $r$ be a positive rational number unequal to 1.  If $c_0, c_1, 
..., c_n, $ are algebraic numbers with $c_0\neq 0$,
and $d$ is an integer, then
$$c_0\pi + \sum_{j=1}^n c_j \log \alpha_j + d\log r $$

is a transcendental number and hence non-zero.}

\begin{proof}  Let $S$ be such that $\log \alpha_j$ ($j \in S$)
is a maximal ${\mathbb Q}$-linearly independent subset of the
$\log \alpha_j$'s so that we can write
$$\sum_{j=1}^n c_j \log \alpha_j = \sum_{j\in S} d_j \log \alpha_j.$$

Our linear form can now be re-written as

$$-ic_0 \log (-1) + \sum_{j\in S} d_j \log \alpha_j + d\log r. $$

By Baker's theorem, this is either zero or transcendental.  We need
to show that the former case cannot arise.  This will follow if
we can show that $$\log (-1), \log \alpha_j \quad (j\in S), \quad \log r$$
are linearly independent over $\mathbb Q$.  But this is indeed the case
since

$$b_0 \log (-1) + \sum_{j\in S}b_j \log \alpha_j + b\log r = 0, $$

for integers $b, b_0, b_j$ ($j\in S$) implies that

$$\prod_{j\in S} \alpha_j^{2b_j} = r^{-2b}. $$

Since the $\alpha_j$'s are units we see that $r^2=1$
Since $r\neq 1$, we have $b=0$.  

By the multiplicative independence of $\alpha_j, j\in S$, we deduce
that $b_j=0$ for all $j \in S$.  Finally, this forces $b_0=0$
and so the numbers are linearly independent over $\mathbb Q$ as
required.

\end{proof}

 \section{\bf The non-vanishing of $L(1,f)$}

\subsection{\textbf{Proof of Theorem 1}}

\smallskip
\noindent

\begin{proof} We first observe that for $\Re(s) >1$, 

\begin{eqnarray*}
L(s,f)&=&\sum_{n=1}^\infty \frac{f(n)}{n^s} \\
&=&\sum_{n=1}^\infty \frac{1}{n^s}\sum_{b=1}^{q-1} \hat{f}(b)\zeta^{bn}_q \\
&=&\sum_{b=1}^{q-1}\hat{f}(b) \sum_{n=1}^\infty \frac{\zeta^{bn}_q}{n^s}, \\
\end{eqnarray*}

keeping in mind that $\hat{f}(q)=0$.  

Now let $s \to 1^+$ to deduce

$$L(1,f) = -\sum_{b=1}^{q-1}\hat{f}(b)\log(1-\zeta^b_q).$$

Using \eqref{2}, we get

\begin{eqnarray*}
L(1,f)&=&-\sum_{b=1}^{q-1}\left(\sum_{j=1}^{\varphi(q)} c_{bj}\omega_j\right)\log(1-\zeta^b_q)\\
&=&-\sum_{j=1}^{\varphi(q)}\omega_j\sum_{b=1}^{q-1}c_{bj}\log(1-\zeta^b_q) \\
&=&\sum_{j=1}^{\varphi(q)}\omega_j\log\alpha_j \quad {\rm (say)}
\end{eqnarray*} 

where $\alpha_j=\prod_{b=1}^{q-1}(1-\zeta^b_q)^{-c_{bj}}\ne 0$ as none of the factors are zero.

Now $L(1,f)=0$ if and only if,

\begin{equation*}
-\omega_1\log\alpha_1=\sum_{j=2}^{\varphi(q)}\omega_j\log\alpha_j
\end{equation*}

if and only if,

\begin{equation*}
\alpha_1^{-1}=\prod_{j=2}^{\varphi(q)}\alpha_j^{\omega_j}   
\end{equation*}

as $\omega_1=1$. Notice that, if some $\alpha_j$ is not equal to 1, then by Baker's theory the 
right hand side of the above identity is a transcendental number. But the left hand side
is a non-zero algebraic number. This contradiction shows that, $L(1,f)=0$ if and only if $\alpha_j=1$ 
for all $j=1,\cdots,\varphi(q)$, i.e.

\begin{equation*}
\underset{{b=1}}{\prod^{q-1}}(1-\zeta^b_q)^{c_{bj}}=1
\end{equation*}

for all $j=1,\cdots,\varphi(q)$. This completes the proof.

\end{proof}

Let us pause to highlight the significance of the previous theorem.
Given a rational-valued function $f$, we can define $\varphi(q)$ functions
$f_j(b)= c_{bj}$ which are all rational-valued.  Even if $f$ is not
rational-valued, a similar analysis leads again to the study of
rational-valued functions and identity \eqref{3}.    

\medskip

For $1\le a <q$ with $(a,q)=1$, let us consider the automorphism $\sigma_a$ of $\Q(\zeta_q)$ given 
by  $\sigma_a(\zeta_q)=\zeta_q^a$. Applying $\sigma_a$ to the identity \eqref{3},
we get 

\begin{equation}\label{4}
\underset{{b=1}}{\prod^{q-1}}(1-\zeta^{ab}_q)^{c_{bj}}=1
\end{equation}

for all $j=1,\cdots,\varphi(q)$.

We define a new function $f_a$ for $1\le a <q$ with $(a,q)=1$, by $f_a(b)=f(ab)$. Clearly $f_a$ is also a 
rational-valued periodic function with period $q$. Hence by the above theorem we have $L(1,f_a)=0$ if and only if 

\begin{equation*}
\underset{{b=1}}{\prod^{q-1}}(1-\zeta^{ab}_q)^{c_{bj}}=1
\end{equation*}

for all $j=1,\cdots,\varphi(q)$.

 Hence an immediate corollary of the Theorem 1 is following:
 
 \medskip
 \noindent

 \textbf{Corollary 1.}

{\sl $L(1,f)=0$ if and only if $L(1,f_a)=0$ for any $1\le a <q$ with $(a,q)=1$. }

 \medskip
 
 \subsection{Proof of Theorem 2}

 \bigskip

 \begin{proof}

 We apply Corollary 1 with  $\sigma_{-1}$ (the complex conjugation) to \eqref{3}, and get 
$L(1,f)=0$ if and only if $L(1,f_{-})=0$ where $f_{-}(x)=f(-x)$. Again, $f$ can be 
written as $$ f=f_e + f_o$$ where 

\begin{equation*}
f_e=\frac{f+f_{-}}{2}
\end{equation*}

and

\begin{equation*}
f_o=\frac{f-f_{-}}{2}
\end{equation*} 

are the even and odd part of $f$ respectively. Hence $L(1,f)=0$ if and only if 
$L(1,f_e)=0$ and $L(1,f_o)=0$. 

\end{proof}

\section{\bf Variation of Okada's Criterion}

As noted earlier, Okada \cite{TO} gave a criterion for the non-vanishing of
$L(1,f)$.  We derive several variations of this here using Dirichlet $L$-series.

\par

\medskip
\noindent

\textbf{Proof of Theorem 3}.

Let $M(q)$ be as before the monoid generated by the prime factors of $q$. Then we have

\begin{eqnarray*}
L(s,f)&=&\sum_{n=1}^\infty \frac{f(n)}{n^s} \\
&=&\sum_{\substack{b\in M(q) \\ (a,q)=1}} \frac{f(ab)}{a^s b^s},
\end{eqnarray*}

since any natural number $n$ can be factored uniquely as $n=ab$ with
$(a,q)=1$ and $b\in M(q)$.

 Thus

\begin{eqnarray*}
L(s,f)&=&\sum_{b\in M(q)} \frac{1}{b^s}\sum_{\substack{a=1\\(a,q)=1}}^\infty\frac{f(ab)}{a^s}\\
&=&\sum_{b\in M(q)} \frac{1}{b^s}\sum_{\substack{a=1\\(a,q)=1}}^\infty\frac{f_b(a)}{a^s}. 
\end{eqnarray*}

Observe that $f_b$ is a function supported
on the coprime residue classes mod $q$. Thus, 
we can write (using the inner product on the group of coprime residue
classes (mod $q$)):  

\begin{equation*}
f_b=\sum_{\chi({\rm mod}~q)}(f_b,\chi)\chi
\end{equation*}

where $\chi$ runs over all the Dirichlet characters mod $q$ and 

\begin{equation*}
(f_b,\chi)=\frac{1}{\varphi(q)}\sum_{\substack{a=1\\(a,q)=1}}^q f_b(a)\overline{\chi}(a).
\end{equation*}

Thus, we have 

\begin{eqnarray*}
L(s,f)=\sum_{b\in M(q)} \frac{1}{b^s}\sum_{\chi({\rm mod}~q)}(f_b,\chi)L(s,\chi).
\end{eqnarray*}

Again, note that for the principal Dirichlet character $\chi_0$ modulo $q$, we have

\begin{eqnarray*}
 L(s,\chi_0)&=&\zeta(s)\prod_{p|q}\left(1-\frac{1}{p^s}\right) \\
 &=& \left\{\frac{1}{s-1}+\gamma+O(s-1)\right\}\left\{\sum_{d|q}\frac{\mu(d)}{d}-
 \left(\sum_{d|q}\frac{\mu(d)}{d}\log d\right)(s-1)+\cdots \right\}  \\
&=&\frac{\phi(q)/q}{s-1}-\sum_{d|q}\frac{\mu(d)}{d}\log d+\gamma\frac{\phi(q)}{q}+O(s-1),          
\end{eqnarray*}

since $\zeta(s)=\left\{\frac{1}{s-1}+\gamma+O(s-1)\right\}$ and 
$$\prod_{p|q}\left(1-\frac{1}{p^s}\right)=\sum_{d|q}{\mu(d)\over d^s}.$$
Hence, we get

\begin{eqnarray*}
L(1,f)&=&\sum_{b\in M(q)} \frac{1}{b}\sum_{\substack{\chi({\rm mod}~q) \\ \chi\ne \chi_0}}(f_b,\chi)L(1,\chi) \\
&+& \underset{s\rightarrow 1^+}{\lim}\sum_{b\in M(q)} \frac{1}{b}\left\{1-(s-1)\log b+\cdots\right\}(f_b,\chi_0)L(s,\chi_0) 
\end{eqnarray*}

Again, for any $c$ with $(c,q)=1$, we have

\begin{eqnarray*}
(f_{cb},\chi)&=&\frac{1}{\varphi(q)}\sum_{\substack{a=1\\(a,q)=1}}^q 
f_{cb}(a)\overline{\chi}(a) \\
&=& \frac{1}{\varphi(q)}\sum_{\substack{t=1\\(t,q)=1}}^q f_b(t)\overline{\chi}(c^{-1}t)\\
&=&\chi(c) (f_b,\chi).            
\end{eqnarray*} 

So, we deduce that

\begin{eqnarray}\label{5}
L(s,f_c)=\sum_{b\in M(q)} \frac{1}{b^s}\sum_{\chi({\rm mod}~q)}\chi(c)(f_b,\chi)L(s,\chi).
\end{eqnarray}

Now, consider the sum

\begin{equation*}
 \sum_{\substack{c=1\\(c,q)=1}}^q L(s,f_c).
\end{equation*}

Note that, for $\chi\ne \chi_0$, we have 

\begin{equation*}
 \sum_{\substack{c=1\\(c,q)=1}}^q \chi(c)=0
\end{equation*}

and

\begin{equation*}
 \sum_{\substack{c=1\\(c,q)=1}}^q \chi_0(c)=\phi(q).
\end{equation*}

Hence from \eqref{5}, we get

\begin{eqnarray*}
\sum_{\substack{c=1\\(c,q)=1}}^q L(s,f_c)&=&\varphi(q)\sum_{b\in M(q)} \frac{1}{b^s}(f_b,\chi_0)L(s,\chi_0) \\
&=&\varphi(q)\sum_{b\in M(q)} \frac{1}{b}\left\{1-(s-1)\log b+\cdots\right\}(f_b,\chi_0)L(s,\chi_0)
\end{eqnarray*}
so that, we have

\begin{eqnarray*}
\sum_{\substack{c=1\\(c,q)=1}}^q L(s,f_c)
&=&\varphi(q)\left[\left(\frac{\varphi(q)}{q}\sum_{b\in M(q)} \frac{(f_b,\chi_0)}{b}\right)\frac{1}{s-1}\right.\\
&+&\sum_{b\in M(q)} \frac{(f_b,\chi_0)}{b}\left(-\sum_{d|q}\frac{\mu(d)}{d}\log d+
\left. \gamma\frac{\varphi(q)}{q}-\frac{\varphi(q)}{q}\log b\right)+O(s-1)\right].
\end{eqnarray*}

Note that the left hand side does not have a pole at $s=1$ and hence from the above equation, we get

\begin{equation}\label{6}
\sum_{b\in M(q)} \frac{(f_b,\chi_0)}{b}=0. 
\end{equation}

Now if $L(1,f)=0$, then $L(1,f_c)=0$ for all $1\le c <q$ with $(c,q)=1$
by Corollary 1. Hence from the above equation, we get 

\begin{equation}\label{7}
 \sum_{b\in M(q)} \frac{(f_b,\chi_0)}{b}\left(-\sum_{d|q}\frac{\mu(d)}{d}\log d+
\gamma\frac{\varphi(q)}{q}-\frac{\varphi(q)}{q}\log b\right)=0.
\end{equation}

Thus using \eqref{6}, the above identity becomes 

\begin{equation}\label{8}
\sum_{b\in M(q)} \frac{(f_b,\chi_0)}{b}\log b=0. 
\end{equation}

Again for any Dirichlet character $\psi\ne \chi_0$ mod $q$, we have

\begin{eqnarray*}
\sum_{\substack{c=1\\(c,q)=1}}^q\overline{\psi}(c)L(s,f_c)&=&\sum_{\substack{c=1\\(c,q)=1}}^q
\overline{\psi}(c)\sum_{b\in M(q)} \frac{1}{b^s}\sum_{\chi(\rm{mod}~q)}\chi(c)(f_b,\chi)L(s,\chi) \\
&=&\sum_{b\in M(q)} \frac{1}{b^s}\sum_{\chi(\rm{mod}~q)}(f_b,\chi)L(s,\chi) 
\sum_{\substack{c=1\\(c,q)=1}}^q\overline{\psi}(c)\chi(c). 
\end{eqnarray*}

Finally, using orthogonality of characters, we get
\begin{eqnarray*}
 \sum_{\substack{c=1\\(c,q)=1}}^q\overline{\psi}(c)L(s,f_c)=\varphi(q)\sum_{b\in M(q)} \frac{(f_b,\psi)}{b^s}L(s,\psi).
\end{eqnarray*}

Now, as before, if $L(1,f)=0$ then $L(1,f_c)=0$ for all $c$ with $(c,q)=1$
by Corollary 1. Hence at $s=1$, left hand side of 
the above identity is 0, so that we have

\begin{equation*}
 \sum_{b\in M(q)} \frac{(f_b,\psi)}{b}L(1,\psi)=0
\end{equation*}

and hence 

\begin{equation}\label{9}
 \sum_{b\in M(q)} \frac{(f_b,\psi)}{b}=0,
\end{equation}

since $L(1,\psi)\ne 0$ by a celebrated theorem of Dirichlet.

Thus, from \eqref{6} and \eqref{9}, we get for all Dirichlet characters $\chi$ mod $q$

\begin{equation*}
 \sum_{b\in M(q)} \frac{(f_b,\chi)}{b}=0.
\end{equation*}

Now multiplying the above identity by $\chi(a)$ and summing over all the Dirichlet 
characters $\chi$ and using the orthogonality of characters, we get 

\begin{equation}\label{10}
 \sum_{b\in M(q)} \frac{f(ab)}{b}=0
\end{equation}

for any $a$ with $(a,q)=1$. Thus \eqref{8} and \eqref{10} prove the only if statement of our theorem.

The if part is also clear from \eqref{6} and the immediate predecessor equation of \eqref{6}. 
This completes the proof of the theorem.

\section{\bf Equivalence of Okada's Criterion}

\bigskip

Theorem 3 looks different from Okada's criterion \cite{TO}.  
 In this section, we prove Theorem 4, which 
shows the equivalence of Theorem 3 and  Okada's criterion. 
\par

Let us begin with the following lemma,
which appears in Okada's paper \cite{TO} without proof. Since the proof is important and relevant 
to our discussion, we give it below.

\medskip

\noindent

\textbf{Lemma 1}.{\sl  The sum

\begin{equation*}
\sum_{j=1}^\infty\sum_{\substack{t=1 \\ p^jt\equiv r({\rm mod}~q)}}^q \frac{1}{p^j}
\end{equation*}
is equal to $\epsilon(r,p)$, where

$$
 \epsilon(r,p)=\begin{cases} 
v_p(r)  & \text{ if $v_p(r)<v_p(q)$,}\\
v_p(q)+\frac{1}{p-1}  & \text{ otherwise. }
\end{cases}
$$}

\noindent

\begin{proof}

 We consider two cases. \\

 \smallskip

 {\bf Case 1:}  $v_p(r)<v_p(q)$. In this case, the congruence $p^jt\equiv r~({\rm mod}~q)$ implies 
 $j\le v_p(r)$, for otherwise, there are no solutions. Thus for  $j\le v_p(r)$, the congruence reduces to 
 $t\equiv p^{v_p(r)-j} r_1~({\rm mod}~p^{v_p(q)-j}q_1)$ where we have written $r=p^{v_p(r)}r_1$, 
 $q=p^{v_p(q)}q_1$ with $(r_1,p)=(q_1,p)=1$. Thus $t$ is uniquely determined $({\rm mod}~p^{v_p(q)-j}q_1)$ 
 for fixed $j$. Now $t$ has precisely $p^j$ lifts $({\rm mod}~q)$. Hence our sum is 

 \begin{equation*}
 \sum_{j=1}^{v_p(r)}\frac{1}{p^j}.p^j=v_p(r)=\epsilon(r,p)
 \end{equation*}
 as desired.

\noindent
\smallskip

{\bf Case 2:}  $v_p(r)\ge v_p(q)$. In this case, if $j\le v_p(q)$, then the congruence $p^jt\equiv r~({\rm mod}~q)$ 
reduces to $t\equiv p^{v_p(r)-j} r_1~({\rm mod}~p^{v_p(q)-j}q_1)$ so that again, the number of lifts of 
$t~({\rm mod}~q)$ is $p^j$. Thus 

\begin{equation*}
 \sum_{j=1}^{v_p(q)}\frac{1}{p^j}.p^j=v_p(q).
 \end{equation*}

Now if $j>v_p(q)$, then our congruence becomes 
$$p^{j-v_p(q)}t\equiv p^{v_p(r)-v_p(q)} r_1~({\rm mod}~q_1). $$

As $(q_1,p)=1$, $t$ is uniquely determined $({\rm mod}~q_1)$. Thus, $t$ has precisely $p^{v_p(q)}$ lifts 
$({\rm mod}~q)$, so that

\begin{equation*}
\sum_{j>v_p(q)} \frac{1}{p^j}.p^{v_p(q)}=\sum_{k=1}^\infty\frac{1}{p^k}=\frac{1}{p-1}.
\end{equation*}

Hence, in this case our sum is

\begin{equation*}
 v_p(q)+\frac{1}{p-1}=\epsilon(r,p).
 \end{equation*}

 This completes the proof.

\end{proof}

\subsection{Proof of Theorem 4}

\medskip

\begin{proof}

 We have

 \begin{equation*}
0=\sum_{b\in M(q)} \frac{(f_b,\chi_0)}{b}\log b. 
\end{equation*}

Let $\Lambda$ be the von Mangoldt function
so that

$$\log b = \sum_{d|b} \Lambda(d). $$

 Then the above identity implies and is implied by the following:

\begin{eqnarray*}
 0&=&\sum_{b\in M(q)} \frac{(f_b,\chi_0)}{b}\sum_{d|b}\Lambda(d) \\
 &=& \sum_{d\in M(q)}\Lambda(d)\sum_{\substack{b\in M(q)\\d|b}} \frac{(f_b,\chi_0)}{b}.
  \end{eqnarray*}

The outer sum is over prime powers in $M(q)$, and we can write each
$d=p^\alpha$, with $p$ a prime divisor of $q$.  Since $d|b$, we write $b=db_1$ 

and get
 \begin{equation}\label{11}
  \sum_{p|q}\log p\left(\sum_{\alpha=1}^\infty \sum_{b_1\in M(q)} \frac{(f_{p^\alpha b_1},\chi_0)}{p^\alpha b_1}\right)=0.
 \end{equation}

Let us study for each prime $p|q$, the coefficient of $\log p$ in \eqref{11}:

\begin{equation*}
\sum_{\alpha=1}^\infty \sum_{b_1\in M(q)} \frac{(f_{p^\alpha b_1},\chi_0)}{p^\alpha b_1}.
\end{equation*}

We write each $b_1$ as $p^\beta c$ with $c\in M(q_1)$ where $q_1=q/p^{v_p(q)}$. Then our sum becomes

\begin{equation*}
 \sum_{\alpha=1}^\infty\sum_{\beta=0}^\infty \sum_{c\in M(q_1)} \frac{(f_{p^{\alpha+\beta}c},\chi_0)}{p^{\alpha+\beta}c}.
\end{equation*}

Substituting the value of $(f_{p^\alpha b_1},\chi_0)$ into the above sum,
we get that this is equal to

\begin{equation*}
\sum_{\alpha=1}^\infty\sum_{\beta=0}^\infty \sum_{c\in M(q_1)}\sum_{\substack{a=1\\(a,q)=1}}^q 
\frac{f(p^{\alpha+\beta}ca)}{p^{\alpha+\beta}c}.  
\end{equation*}

We can collect the powers of $p$, observing that for a fixed $j$, the number of solutions of 
$\alpha+\beta=j$ with $\alpha\ge 1$, $\beta\ge 0$ is precisely equal to $j$. Thus, our sum becomes

\begin{equation*}
\sum_{j=1}^\infty\frac{j}{p^j} \sum_{c\in M(q_1)}\frac{1}{c}\sum_{\substack{a=1\\(a,q)=1}}^q f(p^j ca).  
\end{equation*}

Now let $p^j ca\equiv r \pmod{q}$ so that $(r,q)>1$. Then the above sum becomes

\begin{equation}\label{12}
\sum_{\substack{r=1\\(r,q)>1}}^q f(r)\sum_{j=1}^\infty\frac{j}{p^j} \sum_{c\in M(q_1)}\frac{1}{c}
\sum_{\substack{a=1\\(a,q)=1\\ p^j ca\equiv r ({\rm mod}~q)}}^q 1. 
\end{equation}

We analyse the inner congruence $p^j ca\equiv r \pmod{q}$.
As before, we write $r=p^{v_p(r)}r_1$ so that the congruence becomes

\begin{equation}\label{13}
p^j ca\equiv p^{v_p(r)}r_1  ~ ({\rm mod}~p^{v_p(q)}q_1).
\end{equation}

We consider (as before) two cases.

\noindent
\medskip

 {\bf Case 1:}  $v_p(r)<v_p(q)$. As $ca$ is coprime to $p$, the congruence \eqref{13} implies 
 $j=v_p(r)$ is the only solution for $j$. Thus \eqref{13} reduces to 

 \begin{equation}\label{14}
 ca\equiv r_1  ~ ({\rm mod}~p^{v_p(q)-v_p(r)}q_1). 
 \end{equation}

By the Chinese remainder theorem, this is equivalent to the system 

\begin{equation*}
 ca\equiv r_1  ~ ({\rm mod}~p^{v_p(q)-v_p(r)})
\end{equation*}

\begin{equation*}
 ca\equiv r_1  ~ ({\rm mod}~ q_1).
\end{equation*}

For a given $c$, the second congruence has a solution for $a^{-1}$ if and only if $(r_1,q_1)|c$ in which 
case there are $(r_1,q_1)$ solutions (mod $q_1$). The first congruence has a unique solution for $a^{-1}$ 
(mod $p^{v_p(q)-v_p(r)}$) and thus has $p^{v_p(r)}$ solutions (mod $p^{v_p(q)}$). In total we obtain 
that the inner most sum of \eqref{12} is $p^{v_p(r)}(r_1,q_1)$ in the case  $(r_1,q_1)|c$.

Thus, the two innermost sums in \eqref{12} become (on writing
$c=(r_1, q_1)c_1$ with $c_1 \in M(q_1)$), 

\begin{eqnarray*}
 \sum_{\substack{c\in M(q_1)\\(r_1,q_1)|c}} \frac{1}{c}(r_1,q_1)p^{v_p(r)}
&=&p^{v_p(r)}\sum_{\substack{c_1\in M(q_1)}} \frac{1}{c_1}\\
&=& p^{v_p(r)}\prod_{\substack{p_1|q_1\\p_1 ~{\rm prime}}}\left(1+\frac{1}{p_1}+\frac{1}{p_1^2}+\cdots \right) \\
&=&p^{v_p(r)}\prod_{\substack{p_1|q_1\\p_1 ~{\rm prime}}}\left(1-\frac{1}{p_1}\right)^{-1} \\
&=&\frac{q_1}{\varphi(q_1)}p^{v_p(r)}.
\end{eqnarray*}

Thus, in this case, the total contribution for the three innermost sums in \eqref{12} is 

\begin{equation*}
 v_p(r)\frac{q_1}{\varphi(q_1)}. 
\end{equation*}

Hence, the coefficient of $\log p$ in \eqref{11} is a rational number and is equal to 

\begin{equation*}
 \sum_{\substack{r=1\\(r,q)>1}}^q f(r)v_p(r)\frac{q_1}{\varphi(q_1)}
 =\left(\sum_{\substack{r=1\\(r,q)>1}}^q f(r)\epsilon(r,p)\right)\frac{q_1}{\varphi(q_1)}.
\end{equation*}

\noindent
\textbf{Case 2:} $v_p(r)\ge v_p(q)$. In this case, we must have $j\ge v_p(r)$ and \eqref{13} reduces to 

\begin{equation}\label{15}
c(p^{j-v_p(q)}a) \equiv p^{v_p(r)-v_p(q)}r_1 ~ ({\rm mod}~ q_1).
\end{equation}

Again, for a fixed $c$, this congruence has a solution for $a^{-1}$ if and only if $(r_1,q_1)|c$ in which case 
it has $(r_1,q_1)$ solutions (mod $q_1$). These lift to $\frac{\varphi(q)}{\varphi(q_1)}(r_1,q_1)$ 
solutions (mod $q$).

\medskip

Thus, in this case the inner sums of \eqref{12} become

\begin{equation*}
\sum_{j=v_p(q)}^\infty\frac{j}{p^j} \sum_{\substack{c\in M(q_1)\\(r_1,q_1)|c}}\frac{1}{c}\frac{\varphi(q)}{\varphi(q_1)}(r_1,q_1)
= \sum_{j=v_p(q)}^\infty\frac{j}{p^j}\frac{q_1}{\varphi(q_1)}\frac{\varphi(q)}{\varphi(q_1)}.
\end{equation*}

It is easy to check that

\begin{equation*}
\sum_{j=v}^\infty jX^j = \frac{X^v}{1-X}\left(\frac{X}{1-X}+v\right)
\end{equation*}

so that our inner sums become

\begin{equation*}
 \frac{q_1}{\varphi(q_1)}\frac{\varphi(q)}{\varphi(q_1)}\left(\frac{1}{p}\right)^{v_p(q)}
 \frac{1}{1-1/p}\left(\frac{1/p}{1-1/p}+v_p(q)\right)
\end{equation*}

which is equal to

\begin{equation*}
 \frac{q_1}{\varphi(q_1)}\left(\frac{1}{p-1}+v_p(q)\right).
\end{equation*}

Thus, again we get the coefficient of $\log p$ in \eqref{11} is

\begin{equation*}
 \left(\sum_{\substack{r=1\\(r,q)>1}}^q f(r)\epsilon(r,p)\right)\frac{q_1}{\varphi(q_1)}.
\end{equation*}

Hence, in any case, from our original sum \eqref{11} we obtain 

\begin{equation*}
 \sum_{p|q}\log p\left(\sum_{\substack{r=1\\(r,q)>1}}^q f(r)\epsilon(r,p)\right)\frac{q_1}{\varphi(q_1)}=0.
\end{equation*}

Finally using the unique factorization theorem for the natural numbers, we conclude that

\begin{equation*}
\sum_{\substack{r=1\\(r,q)>1}}^q f(r)\epsilon(r,p)=0
\end{equation*}

for every prime divisor $p$ of $q$.

The converse is also clear from the calculation that
led to the penultimate equation. This completes the proof of the theorem.

\end{proof}

\section{\bf The theorem of Baker, Birch and Wirsing revisited}

\bigskip

Finally, in this section we give a new proof of the theorem of Baker, Birch
and Wirsing \cite{BBW}.  

For the proof of Proposition 1, we shall need the following lemma (see S. Lang \cite{SL}, p.548).

\smallskip
\noindent

\textbf{Lemma 3}. {\sl Let $G$ be any finite abelian group of order $n$ and $F:G\rightarrow \C$ be any complex-valued function
on $G$. The determinant of the $n\times n$ matrix given by ($F(xy^{-1})$) as $x, y$ range over the group elements is 
called the {\it Dedekind determinant} and is equal to

\begin{equation*}
 \prod_\chi\left(\sum_{x\in G} \chi(x)F(x)\right),
\end{equation*}

where the product is over all characters $\chi$ of $G$. }

We also need the following lemma.

\smallskip
\noindent

\textbf{Lemma 4}. {\sl Let $\psi$ be the digamma function as mentioned in section 2. The sum
$$\sum_{\substack{a=1\\(a,q)=1}}^q\psi(a/q) < -\gamma \varphi(q),$$ 

and is hence non-zero.}

\smallskip

\noindent
\begin{proof}

 First notice that, for any $1\le a<q$ we have
 \begin{equation*}
 \sum_{n=1}^\infty{1\over qn(qn+a)}<{1\over a^2} \sum_{n=1}^\infty{1\over n(n+1)}={1\over a^2}
 \end{equation*}

as $qn(qn+a)>an(an+a)$ and the last series of the right hand side of the above inequality is telescopic.  Hence, we get

\begin{equation}\label{16}
 \sum_{n=1}^\infty{1\over n(qn+a)}<{q\over a^2}. 
 \end{equation}

Again we know that ( see \eqref{di})

$$-\psi(x)=\gamma+{1\over x}+\sum_{n=1}^\infty\left({1\over {n+x}}-{1\over n}\right).$$

Hence, we have

$$-\psi(a/q)=\gamma+{q\over a}+\sum_{n=1}^\infty\left({1\over {n+a/q}}-{1\over n}\right)$$

and so that,

$$-\psi(a/q)=\gamma+{q\over a}-\sum_{n=1}^\infty{a\over {n(qn+a)}}.$$

Thus, using \eqref{16} we get $\psi(a/q)<-\gamma$ for any $1\le a<q$. Hence 

$$\sum_{\substack{a=1\\(a,q)=1}}^q\psi(a/q)<-\gamma \varphi(q).$$ This completes the proof. 

\end{proof}

\subsection{\textbf{Proof of Proposition 1}}

We first suppose that $f(q)=0$ and indicate later how this condition
can be removed.

\begin{proof}

 Let $\F=\Q(f(1),\cdots,f(q))$. Then condition $(ii)$ implies $\F\cap\Q(\zeta_q)=\Q$ and hence 
 $[\F(\zeta_q):\F]=\varphi(q)$.  

Now if $L(1,f)=0$, then applying $\sigma_c \in Gal(\F(\zeta_q)/\F)$ and using similar arguments as 
in Theorem 1, we get

 \begin{equation*}
  L(1,f_c)=0
 \end{equation*}

for all $1\le c < q$ with $(c,q)=1$. Here $\sigma_c$ is defined by the rule $\sigma_c(\zeta_q)=\zeta_q^c$

and $f_c(a)=f(ca)$.

Again we know that (see Theorem 16 of \cite{MS1})

\begin{equation*}
L(1,f)=-\frac{1}{q}\sum_{a=1}^q f(a)\psi(a/q). 
\end{equation*}

Thus, we have

\begin{equation*}
 \sum_{a=1}^q f(ca)\psi(a/q)=0
\end{equation*}

for all $1\le c < q$ with $(c,q)=1$.

Rewriting the above equation, we get

\begin{equation*}
\sum_{\substack{b=1 \\ (b,q)=1}}^q f(b)\psi\left(\frac{c^{-1}b}{q}\right)=0 
\end{equation*}

for all $1\le c < q$ with $(c,q)=1$.

Thus we get a matrix equation with $M$ being the $\varphi(q)\times\varphi(q)$ matrix whose $(b,c)$-th entry is given
by $\psi\left(\frac{c^{-1}b}{q}\right)$. Then by the evaluation of the Dedekind determinant as in the Lemma 3,
we get

\begin{equation*}
{\rm det}(M)=\left( \sum_{\substack{a=1\\(a,q)=1}}^q\psi(a/q)\right)
\prod_{\chi\ne \chi_0} q L(1,\chi). 
\end{equation*}

Now using the Lemma 4 and noting that $L(1,\chi)\ne 0$, we get 
${\rm det}(M)\ne 0$.

Thus the matrix $M$ is invertible and hence we have 

\begin{equation*}
f(b)=0, ~1\leq b< q, ~(b,q)=1. 
\end{equation*}

Hence $f$ is identically zero, which is a contradiction to the hypothesis of the theorem.  Finally, 
to treat the case that $f(q) \neq 0$, we adopt a
simple method of Murty and Murty (see section 6 of \cite{murty-murty}).

 This technique
first appears in the prime case in Murty \cite{RM}.  Without loss of
generality, we may suppose that $f$ is integer-valued.  
We define a function $g$ (mod $q$) such that $g(a)=1$ if $(a,q)=1$
and $g(q)=-\varphi(q)$ so that

$$qL(1,g) = -\gamma \varphi(q) - \sum_{(a,q)=1} \psi(a/q) > 0, $$

by Lemma 4.  

As in \cite{murty-murty}, 

$$L(1,g) = {1 \over q} \sum_{d|q} \mu(q/d)d\log d = \sum_{j=1}^\infty S_j, $$
with each $S_j>0$ and $S_j=O(1/j^2)$.  Now define 
$$F(a) = -g(0)f(a) + f(0)g(a), $$

so that $F$ is an integer valued function defined only on the coprime
residue classes (mod $q$).  Thus, $L(1,F) = -g(0)L(1,f) + f(0)L(1,g)$
so that if $L(1,f)=0$, then $L(1,F)=f(0)L(1,g)$.  But $qL(1,g)$ is a ${\mathbb Z}$-linear 
form in logarithms of natural numbers so that  $\exp(qL(1,g))$
is a positive rational number $r$ (say) which is greater than 1.  
On the other hand, $L(1,F) = L(1, F_e) + L(1, F_o)$
where $F_e$ and $F_o$ are the even and odd parts of $F$.  Now, $L(1,F_o)$
is an algebraic multiple of $\pi$, and $L(1,F_e)$ is a linear form
in logarithms of multiplicatively independent real units.  
By Baker's theorem, if $L(1,F)=f(0)L(1,g)$,
we deduce that $\pi$, logarithms of multiplicatively independent real
units and the logarithm of the rational number $r$ are dependent over the
rationals (and hence over the integers).  But this is a contradiction
to Proposition 3.

This completes the proof.

\smallskip
\end{proof}

\bigskip

\noindent
{\bf Acknowledgements.}

We thank Sanoli Gun and Purusottam Rath for their comments on an earlier version of this paper.
We also thank the referee for helpful comments that improved the quality of the paper.

\end{document}